\documentclass[12pt,twoside]{amsart}
\usepackage{amssymb}

\nonstopmode

\numberwithin{equation}{section}
\font\tengothic=eufm10 scaled\magstep 1
\font\sevengothic=eufm7 scaled\magstep 1
\newfam\gothicfam
      \textfont\gothicfam=\tengothic
      \scriptfont\gothicfam=\sevengothic

\newtheorem{theorem}{Theorem}[section]
\newtheorem{lemma}[theorem]{Lemma}
\newtheorem{proposition}[theorem]{Proposition}
\newtheorem{corollary}[theorem]{Corollary}

\theoremstyle{definition}
\newtheorem{definition}[theorem]{Definition} 
\newtheorem{remark}[theorem]{Remark}
\newtheorem{example}[theorem]{Example}

\newtheorem{notation}[theorem]{Notation}
\newtheorem{question}[theorem]{Question}
\newtheorem{construction}[theorem]{Construction}

\newcommand\rank{\operatorname{rank}}

\newcommand\depth{\operatorname{depth}}

\newcommand\im{\operatorname{im}}

\newcommand{\TT}{\operatorname{Tor}}
\newcommand{\TR}{\operatorname{Tor}^R}

\newcommand{\proj}[1]
{ \mathchoice
           { {\mathbb P}^{#1} }
           { {\mathbb P}^{#1} }
           { {\mathbb P}^{#1} }
           { {\mathbb P}^{#1} }
         }

\newcommand{\s}{\; | \;}
\newcommand{\mif}{\mbox{if} ~}

\newcommand{\timesl}{\hbox{\scriptsize $(\times L)$}}

\newcommand{\fm}{{\mathfrak m}}

\newcommand{\WLP}{Weak Lefschetz property}
\newcommand{\tor}{\operatorname{tor}}
\newcommand{\tR}{\operatorname{tor}^R}
\newcommand{\fa}{{\mathfrak a}}
\newcommand{\fb}{{\mathfrak b}}
\newcommand{\Soc}{\operatorname{Soc}}

\newcommand {\ZZ}{\mathbb{Z}}

\newcommand{\kindalong}{\ \hbox{$\hbox to .35in{\rightarrowfill}$} \ }
\newcommand{\reallylong}{\hbox{$\hbox to .5in{\rightarrowfill}$}  }
\def\ra{\rightarrow}
\def\X{{\mathbb X}}
\def\Y{{\mathbb Y}}

\begin{document}
\title[Weak and Strong Lefschetz Properties]{The Weak and Strong Lefschetz
Properties for Artinian $K$-algebras}

\author{Tadahito Harima, Juan C. Migliore,
Uwe Nagel, Junzo Watanabe } 
\address{Department of Information Science, 
Shikoku University, 
Tokushima 771-1192, Japan} 
\email{harima@keiei.shikoku-u.ac.jp} 
\address{Department of Mathematics,
        University of Notre Dame,
        Notre Dame, IN 46556,
        USA}
\email{migliore.1@nd.edu} 
\address{Fachbereich Mathematik und Informatik, Universit\"at 
Paderborn, D--33095 Paderborn, Germany}
\email{uwen@uni-paderborn.de} 
\address{Department of Mathematical Sciences, Tokai University, Hiratsuka
  259-1292, Japan}
\email{junzowat@ss.u-tokai.ac.jp} 


\begin{abstract}
Let $A = \bigoplus_{i \geq 0} A_i$ be a standard graded Artinian
$K$-algebra, where $\hbox{char } K = 0$.  Then $A$ has the Weak Lefschetz
property if there is an element $\ell$ of degree 1 such that the
multiplication $\times \ell : A_i \rightarrow A_{i+1}$ has maximal rank, for
every $i$, and $A$ has the Strong Lefschetz property if  $\times \ell^d :
A_i \rightarrow A_{i+d}$ has maximal rank for every $i$ and $d$.

The main results obtained in this paper are the following.

\begin{enumerate}
\item {\em Every} height three complete intersection has the Weak
Lefschetz property.  (Our method, surprisingly, uses rank two vector
bundles on $\proj{2}$ and the Grauert-M\"ulich theorem.)

\item We give a complete characterization (including a concrete
construction) of the Hilbert functions that can occur for $K$-algebras
with the Weak or Strong Lefschetz property (and the characterization is
the same one!).

\item We give a sharp bound on the graded Betti numbers (achieved by our
construction) of Artinian $K$-algebras with the Weak or Strong Lefschetz
property and fixed Hilbert function.  This bound is again the same for
both properties!  Some Hilbert functions in fact {\em force} the algebra to
have the maximal Betti numbers.

\item {\em Every} Artinian ideal in $K[x,y]$ possesses the Strong
Lefschetz property.  This is false in higher codimension.

\end{enumerate}
\end{abstract}


\maketitle

\section{Introduction}

Let $A$ be a graded Artinian algebra over some field $K$ (which we will
restrict shortly). Then $A$ has the {\em Weak Lefschetz property} (sometimes
called the Weak Stanley property) if there is an element $\ell$ of degree
1 such that the multiplication $\times \ell : A_i \rightarrow A_{i+1}$ has
maximal rank, for every $i$.  We say that $A$ has the {\em Strong Lefschetz
property} if there is an element $\ell$ of degree 1 such that the
multiplication $\times \ell^d : A_i \rightarrow A_{i+d}$ has maximal rank
for every $i$ and $d$.  If $A = R/I$, where $R$ is a polynomial ring and
$I$ is a homogeneous ideal, then sometimes we will abuse notation and refer
to the Weak or Strong Lefschetz properties for $I$ rather than for $A$.
These are both fundamental properties and have been studied by many authors,
especially when $A$ is Gorenstein (e.g.\ \cite{boij}, \cite{harima},
\cite{iarrobino}, \cite{ikeda}, \cite{MN3}, \cite{stanley1},
\cite{watanabe}, \cite{watanabe1}, \cite{watanabe2}).

Throughout this paper, unless specified otherwise, we assume that we
work over a field of characteristic zero.  This paper began with a study of
the Weak Lefschetz property for complete intersections of height three, and
grew to a study of Artinian ideals of arbitrary codimension.  Our original
interest in the subject was to try to get a handle on ``how many'' Artinian
complete intersections possess this natural property.  However, a further
motivation comes from the fact that this property can be translated into
(at least) two other natural questions.

First, suppose  that $F_1, F_2, \dots, F_n$ is a homogeneous
complete intersection in the $n$ dimensional polynomial ring $R$.   Then the
minimal free resolution of the ideal $(F_1, \dots, F_n)$ is well
understood; namely it is obtained as the Koszul complex.  However, the
graded Betti numbers of the minimal free resolution of the ideal $(F_1,
\dots, F_n, L)$, where L is a generic linear form, does not seem to be well
understood. For example, should they be all the same, depending only on the
degrees of the generators and not on the generators themselves, as long as
they are a regular sequence of given degrees plus a generic element? ( We
could also ask the same question for $L^d$ in the place for $L$.)  The
connection between the Weak Lefschetz property and this question is
discussed in the last part of section 2, and we give a complete answer
(Corollary \ref{samerestr}) when $n=3$.

One other problem concerns the generic initial ideal, $\hbox{gin}(I)$, of a
complete  intersection $I$, i.e.\ the initial ideal of $I$ with respect to
generic variables (cf.\ for instance \cite{eisenbud}).   It is well known
that $\hbox{gin}(I)$ is Borel-fixed.   But if $I$ is a complete
intersection and if we fix a monomial order, is the Borel-fixed ideal
$\hbox{gin}(I)$ unique? Or are there two complete intersections $I$ and
$J$  such that $\hbox{gin}(I)$ and
$\hbox{gin}(J)$ are different Borel-fixed ideals with the same Hilbert
function?
These questions seem to be open since if $\hbox{gin}(I)$ is unique with
respect to the reverse lex order then it would imply the Strong Lefschetz
property of all complete intersections of those degrees. Since a
Borel-fixed ideal is unique in codimension two (for a fixed Hilbert
function) the Strong Lefschetz property can be proved in this case
(Proposition \ref{SLP for codim 2}).

It should also be mentioned that Stanley and others have made deep
connections between the Weak and Strong Lefschetz properties and questions
in  combinatorics \cite{stanley1}, \cite{stanley2}.  For example, the Weak
Lefschetz property was the crucial ingredient in Stanley's part of the
characterization of the $f$-vectors of simplicial polytopes.  Thus, we are
exploring in this paper also the restrictions on the possible Hilbert
functions and graded Betti numbers imposed by the presence of the Weak or
Strong Lefschetz property.

It was noticed by Stanley \cite{stanley2} and independently by the fourth
author \cite{watanabe1} that any monomial complete intersection (in any
number of variables) has the Strong Lefschetz property, and the
fourth author proved that in any codimension, ``most'' Artinian Gorenstein
rings with fixed socle degree possess the Strong  Lefschetz property
(\cite{watanabe1}, Example 3.9).  We remark (following \cite{iarrobino})
that Stanley's proof used the idea of recognizing $A = R/I$ as the
cohomology ring of a product $X$ of projective spaces, and then using the
hard Lefschetz theorem for the algebraic variety $X$.  The fourth author
noticed that if follows from the representaion theory of the Lie algebra
$sl(2)$.

Yet even in codimension 3, we do not have a clear idea of which Artin
Gorenstein rings possess this property, and in particular whether all of
them do.  The (apparently) simplest situation is for height 3 complete
intersections in $R=K[x_1,x_2,x_3]$.  Until now the most general result
for this case is again due to the fourth author.  Suppose that the
generators of the complete intersection $I$ have degrees $2 \leq d_1 \leq
d_2 \leq d_3$.  Then it was proved in \cite{watanabe2} that if $d_3 > d_1
+ d_2 - 2$ then  $R/I$ has the Weak Lefschetz property.  But for arbitrary
complete intersections, even the case of three polynomials of degree 4 had
been open.

The first main result of this paper (Theorem \ref{mainresult}) is that {\em
all} Artinian complete intersections in $K[x_1,x_2,x_3]$ have the Weak
Lefschetz property.   It is a somewhat surprising result.
Indeed, it was known to be a very difficult problem among the experts, and
at times it seemed more natural to seek a counter-example rather than
to try to prove it!  We are able to give a relatively simple proof by
translating the problem to one of vector bundles on $\proj{2}$ and
invoking a deep theorem due to Grauert and M\"ulich.

This part of the paper was inspired by \cite{watanabe2}, but as mentioned
earlier, our techniques are completely different from those of the papers
cited above.  Because we apply the Grauert-M\"ulich theorem, we are forced to
assume characteristic zero (as indeed was done in \cite{watanabe2}).  In
fact, the Weak Lefschetz property does not hold for all complete
intersections in characteristic $p$; see Remark \ref{charp}.

As a further illustration of the power of our approach, we give a simple
proof (Corollary \ref{junzomainresult}) of the main result of
\cite{watanabe2}.

In the third section of the paper we do not assume that $\hbox{char } K =
0$.  We consider graded Artinian
$K$-algebras which are not necessarily complete intersections.  Here we
produce (Construction \ref{constr})  a particular graded Artinian
$K$-algebra, which allows us to give a necessary and sufficient condition
for a sequence of integers to be the Hilbert function of a graded Artinian
$K$-algebra with the Weak Lefschetz property (Proposition \ref{exex}).  We
also answer several natural questions about the minimal free resolutions of
algebras with the Weak Lefschetz property.

Our second main result (Theorem \ref{thm-bounds}) shows that if we fix an
allowable Hilbert function then there is a sharp upper bound on the graded
Betti numbers among $K$-algebras having the Weak Lefschetz property.
Indeed, this bound is achieved by the algebra produced by Construction
\ref{constr}, once we refine the construction slightly.  This result is
analogous to the main result of \cite{MN3}, which proved it for Gorenstein
ideals with the Weak Lefschetz property (see also \cite{GHS1}).  As a
corollary we show that there are Hilbert functions which occur for
$K$-algebras with the Weak Lefschetz property and for which this property
forces the graded Betti numbers to be the maximal ones.

In section 4 we again assume $\hbox{char } K = 0$.  We consider the
Strong Lefschetz property, namely that there exists a linear form $\ell$
such that for each
$d$, the multiplication
$\times \ell^d : A_i \rightarrow A_{i+d}$ has maximal rank, for every $i$.
This condition implies the Weak Lefschetz property, but is not equivalent to
it in general.  We show that these conditions are both automatic in
codimension two, however.

Since there are algebras with the Weak Lefschetz property but
not the Strong Lefschetz property, one might guess that the imposition of
the Strong Lefschetz property reduces the number of possible Hilbert
functions.  However, we are able to show that with another slight refinement
of Construction \ref{constr}, that algebra has the Strong Lefschetz
property.  This yields the surprising result that a Hilbert function occurs
among algebras with the Weak Lefschetz property if and only if it occurs
among algebras with the Strong Lefschetz property.  Furthermore, the
extremal graded Betti numbers for algebras with the Weak Lefschetz property
also occur among algebras with the Strong Lefschetz property.

Our results have some consequences for the punctual Hilbert scheme.
 Since by semicontinuity the Weak and Strong Lefschetz properties
are open properties, it follows that the general point of a  component has
the Strong (resp.\ Weak) Lefschetz property if and only the
component has one point with the Strong (resp.\ Weak) Lefschetz property.
Moreover, we know precisely the possible Hilbert functions of the
$K$-algebras corresponding to such a general point.

The second author would like to thank Chris Peterson and Mohan Kumar for
a helpful conversation about vector bundles.


\section{The Weak Lefschetz Property for height three complete intersections}

Let $R = K[x_1,x_2,x_3]$, where $K$ is a field of characteristic zero.
Initially we will assume that $K$ is algebraically closed, in order to
freely use the results of \cite{OSS}.  However, we note in Corollary
\ref{not alg closed} and beyond that our results hold without
that assumption.

Let $I$ be a complete intersection ideal of $R$
generated by homogeneous elements $F_1,F_2,F_3 \in R$ of degrees
$d_1,d_2,d_3$ respectively, and
$d_1 \leq d_2 \leq d_3$.   The minimal free resolution for $R/I$ has the
form
\begin{equation}\label{RESOL}
\begin{array}{c}
0 \rightarrow R(-d_1-d_2-d_3) \rightarrow  {\mathbb F}_2 \kindalong
{\mathbb F}_1 \stackrel{[F_1,F_2,F_3]}{\reallylong} R \rightarrow R/I
\rightarrow 0 \\
\hskip .15in \searrow \hskip .25in \nearrow \\
\hskip .1in E \\
\hskip .15in \nearrow \hskip .25in \searrow \\
\hskip .1in 0 \hskip .7in 0
\end{array}
\end{equation}
where ${\mathbb F}_2 = R(-d_2-d_3) \oplus R(-d_1-d_3) \oplus R(-d_1-d_2)$
and ${\mathbb F}_1 = R(-d_1) \oplus R(-d_2) \oplus R(-d_3)$.  Sheafifying,
we get the following two exact sequences:
\begin{equation} \label{SES1}
0 \rightarrow {\mathcal E} \rightarrow {\mathcal F}_1
\stackrel{[F_1,F_2,F_3]}{\reallylong} {\mathcal O}_{\proj{2}}
\rightarrow 0
\end{equation}
and
\begin{equation} \label{SES2}
0 \rightarrow {\mathcal O}_{\proj{2}} (-d_1-d_2-d_3) \rightarrow {\mathcal
F}_2 \rightarrow {\mathcal E} \rightarrow 0,
\end{equation}
where $\mathcal E$ is locally free (since $I$ is Artinian) of rank two,
${\mathcal F}_1 = {\mathcal O}_{\proj{2}} (-d_1) \oplus {\mathcal
O}_{\proj{2}} (-d_2) \oplus {\mathcal O}_{\proj{2}} (-d_3)$ and ${\mathcal
F}_2 = {\mathcal O}_{\proj{2}} (-d_2-d_3) \oplus {\mathcal O}_{\proj{2}}
(-d_1-d_3) \oplus {\mathcal O}_{\proj{2}} (-d_1-d_2)$.

We would like a condition which forces $\mathcal{E}$ to be semistable.  We
first consider the case where $d_1+d_2+d_3$ is even.  Choose an integer $d$
so that $2d = d_1+d_2+d_3$.  Notice that $c_1({\mathcal E}) = -d_1-d_2-d_3
= -2d$, so the normalized bundle $\mathcal{E}_{norm}$ is $\mathcal{E}(d)$
(an easy computation, or see for instance \cite{OSS} page 165).   Twisting
the sequence (\ref{SES2}) by $d-1$ we obtain
\begin{equation}\label{computess1}
0 \rightarrow \mathcal{O}_{\proj{2}} (-d-1) \rightarrow
\begin{array}{c}
\mathcal{O}_{\proj{2}} (-d+d_1-1) \\ \oplus \\ \mathcal{O}_{\proj{2}}
(-d+d_2-1) \\
\oplus \\ \mathcal{O}_{\proj{2}} (-d+d_3-1)
\end{array}
 \rightarrow \mathcal{E}_{norm} (-1) \rightarrow 0.
\end{equation}

We now consider the case where $d_1+d_2+d_3$ is odd.  Choose $d$ so that
$2d = d_1+d_2+d_3-1$.  Then again $\mathcal{E}_{norm} = \mathcal{E}(d)$
(again see \cite{OSS} page 165).   Now we have the short exact sequence
\begin{equation} \label{computess2}
0 \rightarrow \mathcal{O}_{\proj{2}} (-d-1) \rightarrow
\begin{array}{c}
\mathcal{O}_{\proj{2}} (-d+d_1-1) \\ \oplus \\ \mathcal{O}_{\proj{2}}
(-d+d_2-1) \\
\oplus \\ \mathcal{O}_{\proj{2}} (-d+d_3-1)
\end{array}
 \rightarrow \mathcal{E}_{norm} \rightarrow 0.
\end{equation}

\newpage

\begin{lemma} \label{semistable}
Let $\mathcal{E}$ be the rank two locally free sheaf obtained above as
the kernel of the map $[F_1,F_2,F_3]$.
\begin{itemize}
\item[1.] Assume $d_1+d_2+d_3$ is even.   If $d_3 < d_1+d_2+2$ then
$\mathcal{E}$ is semistable.

\item[2.] Assume $d_1+d_2+d_3$ is odd.  If $d_3 < d_1+d_2+1$ then
$\mathcal{E}$ is semistable.
\end{itemize}
\end{lemma}

\begin{proof}

When $c_1(\mathcal{E})$ is even and $\mathcal{E}$ has rank two, we know
from \cite{OSS} Lemma 1.2.5 that  $\mathcal E$ is semistable if and only if
$H^0(\proj{n}, \mathcal{E}_{norm}(-1)) = 0$ (since it has rank two).  When
$c_1(\mathcal{E})$ is odd and $\mathcal{E}$ has rank two, stability and
semistability coincide (\cite{OSS} page 166) and the condition for
semistability is $H^0(\proj{2}, \mathcal{E}_{norm}) = 0$.

The two sequences (\ref{computess1}) and (\ref{computess2})  are
exact on global sections.  Hence semistability follows in either case if we
have $-d+d_3-1 < 0$ (where $d$ changes slightly depending on the parity of
$d_1+d_2+d_3$).  The lemma then follows from a simple computation.
\end{proof}

Let $\lambda \cong \proj{1}$ be a general line in $\proj{2}$.  Recall that
every vector bundle on $\proj{1}$ splits, so in particular
$\mathcal{E}|_\lambda \cong \mathcal{O}_{\proj{1}} (e_1) \oplus
\mathcal{O}_{\proj{1}} (e_2)$.  The pair $(e_1,e_2)$ is called the {\em
splitting type} of $\mathcal{E}$.

\begin{corollary}\label{splittingtype}
Let $\mathcal{E}$ be the locally free sheaf obtained above, and assume
that \linebreak
$d_3 < d_1+d_2+1$.  Then the splitting type of $\mathcal{E}$ is
\[
(e_1,e_2) =
\left \{
\begin{array}{ll}
(-d,-d) & \hbox{if } d_1+d_2+d_3 = 2d; \\
(-d, -d-1) & \hbox{if } d_1+d_2+d_3 -1 = 2d.
\end{array}
\right.
\]
\end{corollary}

\begin{proof}
By Lemma \ref{semistable}, $\mathcal{E}$ is semistable.
The theorem of Grauert and M\"ulich (\cite{OSS} page 206, \cite{ein}, page
68) says that in characteristic zero the splitting type of a semistable
normalized 2-bundle $\mathcal{E}_{norm}$ over $\proj{n}$ is
\[
(e_1,e_2) =
\left \{
\begin{array}{ll}
(0,0) & \hbox{if } c_1(\mathcal{E}_{norm}) =0; \\
(0, -1) & \hbox{if } c_1(\mathcal{E}_{norm}) =-1.
\end{array}
\right.
\]
In our case $\mathcal{E}_{norm} = \mathcal{E}(d)$, so a simple
calculation gives the result.
\end{proof}

With this preparation, we now prove the main result of the paper.  We
continue to assume that $K$ is algebraically closed of characteristic zero.

\begin{theorem}  \label{mainresult}
Every height three Artinian complete intersection has the
Weak Lefschetz Property.
\end{theorem}

\begin{proof}
It was shown in \cite{watanabe2} Corollary 3 that if $d_3 \geq d_1+d_2-3$
then $R/I$ has the Weak Lefschetz property.  So without loss of generality
assume that
$d_3 < d_1+d_2-3$.  Note that then Corollary \ref{splittingtype} applies.
To prove the Weak Lefschetz property it is enough to prove injectivity in
the ``first half,'' so we will focus on this.

Let $L$ be a general linear form and let $\bar R = R/L$.  We denote by
$\bar F$ the restriction of a polynomial $F$ to $\bar R$ and by $\bar
{\mathbb F}_1$ the free $\bar R$-module $\bar R(-d_1) \oplus \bar R(-d_2)
\oplus \bar R(-d_3)$.  Consider the multiplication induced by
$L$.  From (\ref{RESOL}) we obtain a commutative diagram

\newcommand{\linmat}
{\scriptsize \left [ \begin{array}{ccc}
L & 0 & 0 \\ 0 & L & 0 \\ 0 & 0 & L
\end{array} \right ]}

\begin{equation} \label{diag1}
\begin{array}{ccccccccccccccc}
&& && 0 && 0 \\
&&&& \downarrow && \downarrow \\
0 & \rightarrow & E(-1) & \rightarrow & {\mathbb F}_1 (-1) &
\stackrel{[F_1 \ F_2 \ F_3]}{\reallylong} & R(-1) &  \rightarrow & R/I (-1)
& \rightarrow & 0\\
&&&& \phantom{M} \downarrow {\scriptstyle M} && \phantom{\timesl}
\downarrow  \timesl &&  \phantom{\timesl} \downarrow  \timesl \\
0 & \rightarrow & E & \rightarrow & {\mathbb F}_1  &
\stackrel{[F_1 \ F_2 \ F_3]}{\reallylong} & R &  \rightarrow & R/I  &
\rightarrow & 0\\
&&&& \downarrow && \downarrow  \\
&&&& \bar {\mathbb F}_1 & \stackrel{[\bar F_1\ \bar F_2 \ \bar
F_3]}{\reallylong} & \bar R & \\
&&&& \downarrow && \downarrow \\
&&&& 0 && 0
\end{array}
\end{equation}
where $M$ is the matrix $\linmat$.  Note that the first vertical exact
sequence is the direct sum of three copies of the exact sequence
\[
0 \rightarrow R(-1) \stackrel{\times L}{\longrightarrow} R \rightarrow
\bar R \rightarrow 0
\]
twisted by $-d_1$, $-d_2$ and $-d_3$ respectively.  The induced map on the
kernels
\[
E(-1) \rightarrow E
\]
is just multiplication by $L$.

Let $\lambda$ be the line in $\proj{2}$ defined by $L$.  Invoking the
Snake Lemma and using the fact that the sheafification of $R/I$ is 0, the
sheafified version of (\ref{diag1}) is

\begin{equation} \label{diag2}
\begin{array}{ccccccccccccccc}
&& 0 && 0 && 0 \\
&& \downarrow && \downarrow && \downarrow \\
0 & \rightarrow & {\mathcal E}(-1) & \rightarrow & {\mathcal F}_1 (-1) &
\stackrel{[F_1 \ F_2 \ F_3]}{\reallylong} & {\mathcal O}_{\proj{2}} (-1) &
\rightarrow &  0\\
&& \phantom{\timesl} \downarrow \timesl && \phantom{M} \downarrow
{\scriptstyle M} && \phantom{\timesl} \downarrow  \timesl  \\
0 & \rightarrow & {\mathcal E} & \rightarrow & {\mathcal F}_1  &
\stackrel{[F_1 \ F_2 \ F_3]}{\reallylong} & {\mathcal O}_{\proj{2}} &
\rightarrow & 0\\
&&\downarrow && \downarrow && \downarrow  \\
0 & \rightarrow & {\mathcal E}|_{\lambda} & \rightarrow & \bar {\mathcal
F}_1 &
\stackrel{[\bar F_1\
\bar F_2 \ \bar F_3]}{\reallylong} & {\mathcal O}_\lambda & \rightarrow & 0
\\ && \downarrow && \downarrow && \downarrow \\
&& 0 && 0 && 0
\end{array}
\end{equation}
By Corollary \ref{splittingtype},
\[
{\mathcal E}|_\lambda \cong \left \{
\begin{array}{rl}
{\mathcal O}_\lambda (-d)^2, & \hbox{ if $d_1+d_2+d_3 = 2d$;} \\ \\
{\mathcal O}_\lambda (-d) \oplus {\mathcal O}_\lambda (-d-1), & \hbox{ if
$d_1+d_2+d_3-1 = 2d$.}
\end{array}
\right.
\]

Let $\bar I$ be the ideal $(\bar F_1, \bar F_2, \bar F_3)$ in $\bar R$.
Taking global sections on the last line of (\ref{diag2}) gives
\[
0 \rightarrow \bigoplus_{i=1}^2 \bar R(-e_i) \rightarrow \bar {\mathbb F}_1
\rightarrow \bar I \rightarrow 0
\]
where $|e_1-e_2| = 0$ or 1 according to whether $d_1+d_2+d_3$ is even or
odd, respectively.  It was observed in \cite{watanabe2} (Remark on page
3165) that this implies that $R/I$ has the Weak Lefschetz property.
However, for completeness we will sketch the argument.  We will treat only
the case
$d_1+d_2+d_3$ even, leaving the other case to the reader.

We have the exact sequence
\begin{equation} \label{restrseq}
0 \rightarrow \bar R(-d)^2 \rightarrow
\begin{array}{c}
\bar R(-d_1) \\ \oplus \\ \bar R(-d_2) \\ \oplus \\ \bar R(-d_3)
\end{array}
\stackrel{[\bar F_1\ \bar F_2 \ \bar F_3]}{\reallylong} \bar I \rightarrow
0
\end{equation}
where $d = \frac{d_1+d_2+d_3}{2}$.

As noted earlier, we only have to show that multiplication by $L$ is
injective on the ``first half"  of $R/I$.  The socle degree of $R/I$ is
$d_1+d_2+d_3-3$, so we have to show that the multiplication
\[
(R/I)_j \stackrel{\times L}{\longrightarrow} (R/I)_{j+1}
\]
is injective for $j \leq \frac{d_1+d_2+d_3}{2} -2 = d-2$.  We will show it
to be true for $j = d-2$, and from the form of the proof it will be clear
that it holds also for smaller $j$.

The kernel of $(\times L)$ is $[I:_R L]$, so if $(\times L)$ is not
injective we have an element $F \in R_{d-2}$, $F \notin I$, such that $LF
\in I_{d-1}$.  That is, we have forms $A_i$, $1 \leq i \leq 3$, with
$\deg A_i = d-1 -d_i$ and
\[
LF - A_1F_1 - A_2F_2 - A_3F_3 = 0.
\]
Restricting this syzygy to $\bar R$ gives
\[
\bar A_1 \bar F_1 + \bar A_2 \bar F_2 + \bar A_3 \bar F_3 = 0.
\]
But (\ref{restrseq}) says that the smallest possible syzygies come from
polynomials of degree $d-d_i$,  $1 \leq i \leq 3$, so this is a
contradiction.  As noted, this works equally well to prove injectivity for
all  $j \leq d-2$.
\end{proof}

\begin{corollary} \label{not alg closed}
Let $K$ be a field of characteristic zero which is not necessarily
algebraically closed.  Then every height three Artinian complete
intersection in $K[x_1,x_2,x_3]$ has the Weak Lefschetz property.
\end{corollary}

\begin{proof}
The Weak Lefschetz property for a graded Artinian $K$-algebra $A$ is
equivalent to the statement that for a general linear form $\ell$, the
Hilbert function of $A/\ell A$ is just the positive part of the first
difference of the Hilbert function of $A$.  But this does not change under
extension of the base field, so the result follows from Theorem
\ref{mainresult}.
\end{proof}

Using the same methods, we can also give a new proof of the main
result of \cite{watanabe2}.  As above, we can assume that $K$ is
algebraically closed initially, but the rest of the results of this section
do not need this assumption.

\newpage

\begin{corollary} \label{junzomainresult}
Let $R=K[x_1,x_2,x_3]$, $I = (F_1,F_2,F_3)$ a complete intersection in $R$,
$d_i = \deg F_i$ for $i=1,2,3$, $L$ a general linear form, $\bar R = R/LR$
and $\bar I = (I+LR)/LR$.  Then the following are equivalent:
\begin{itemize}
\item[(i)] $\mu (\bar I) = 3$, where $\mu$ is the minimal number of
generators;
\item[(ii)] $d_3 \leq d_1 + d_2 -2$.
\end{itemize}
\end{corollary}

\begin{proof}
For completeness we repeat the proof from \cite{watanabe2} of the fact that
(i) implies (ii).  Since $L$ is general, $F_1,F_2$ and $L$ are a regular
sequence, and the socle degree of $R/(F_1,F_2,L)$ is $d_1+d_2-2$.  If (ii)
is not true then $F_3$ is contained in the ideal $(F_1,F_2,L)$, so $\bar
F_3$ is contained in $((F_1,F_2) +LR)/LR$, contradicting (i).

The hard part of the proof is the converse, which we prove using our
approach.  We have from Corollary \ref{splittingtype} that the splitting
type of $\mathcal E$ is
\[
(e_1,e_2) =
\left \{
\begin{array}{ll}
(-d,-d) & \hbox{if } d_1+d_2+d_3 = 2d; \\
(-d, -d-1) & \hbox{if } d_1+d_2+d_3 -1 = 2d.
\end{array}
\right.
\]
With this definition of $d$, a simple calculation  gives that
\[
\begin{array}{lrcl}
\hbox{If } d \hbox{ is even  then} & d_3 < d & \Leftrightarrow & d_3 < d_1
+ d_2 ; \\
\hbox{If } d \hbox{ is odd  then} & d_3 < d & \Leftrightarrow & d_3 < d_1
+ d_2 -1.
\end{array}
\]
So in either case, if (ii) holds then $d_3 < d$.  But the splitting type
gives exactly the leftmost free module in the short exact sequence
(\ref{restrseq}), and the fact that $d_3 < d$ means that no splitting can
occur in the resolution.
\end{proof}

We now apply these ideas to the question of minimal free resolutions.  In
particular, suppose $I = (F_1,F_2,F_3)$ is a complete intersection in $R =
K[x_1,x_2,x_3]$ and $F$ is a general form of degree $d$.  What can be the
possible minimal free resolutions of the ideal $(I,F)$?  Does it depend
only on the degrees of the generators of $I$, or does the choice of the
complete intersection itself play a role?  We can answer this question
when $F$ has degree 1, which in any case was an open question.  To be
consistent with notation, we write $L$ for this general linear
polynomial.  We begin with a lemma.

\begin{lemma} \label{ht2}
Let $I \subset R = K[x_1,x_2,x_3]$ be an Artinian ideal.  Then there
exists a Cohen-Macaulay height two ideal $J \subset R$ such that $J+(L) =
I+(L)$.  $J$ can even be taken to be reduced.
\end{lemma}

\begin{proof}
Let $I = (F_1,\ldots,F_k)$.
We know that $I+(L)/(L) = (\bar F_1,\ldots,\bar F_k)$ is Artinian in
$\bar R = R/(L)$, hence Cohen-Macaulay of height 2. After a change of
coordinates, we can assume that $L = x_3$, hence we obtain polynomials
$G_1,\ldots,G_k \in K[x_1,x_2]$ by canceling all monomials in $F_1,\ldots,F_k$
which are a multiple of $x_3$.  Then viewing these polynomials in $R$
gives the first result.  This ideal is not reduced, however.  But it has a
Hilbert-Burch matrix, whose entries are all polynomials in $x_1,x_2$.
Using standard lifting techniques one can obtain a reduced scheme $J$ with
the desired property.  (A more geometric use of this trick may be found in
\cite{chiantini-orecchia}.)
\end{proof}

Note that the preceding lemma trivially implies that all the graded Betti
numbers (over $R/(L)$) of the reduction of $J$ modulo $L$ are the same as
those of the reduction of $I$ modulo $L$.  However, in general we are not
able to say what these Betti numbers are, or what the Betti numbers of the
ideal $I+(L)$ are (over $R$), or even what the Hilbert function is.
Nevertheless, in the case of complete intersections we can say something
much stronger, thanks to our results above.

\begin{corollary} \label{samerestr}
Let $I = (F_1,F_2,F_3) \subset R$ be a complete intersection.  Then there
is a (reduced) arithmetically Cohen-Macaulay ideal $J = (G_1,G_2,G_3)
\subset R$ such that $\deg G_i = \deg F_i = d_i$ for $i=1,2,3$ and such
that $J+(L) = I+(L)$.  Furthermore:
\begin{itemize}
\item[a.] If $d_3 \leq d_1 + d_2 -2$ then $J$ is an
almost complete intersection with minimal generators given by the $G_i$.
Let $d$ be defined by
\[
\left \{
\begin{array}{ll}
d_1 + d_2+d_3 = 2d & \hbox{ if $d_1+d_2+d_3$ is even} \\
d_1 + d_2+d_3 -1 = 2d & \hbox{ if $d_1+d_2+d_3$ is odd}
\end{array}
\right.
\]
If $d_1+d_2+d_3$ is even then the minimal free resolution of $R/(I+(L))$ is
given by
\[
\begin{array}{cccccccccccccccccccccccc}
&&&&R(-d_1-1) && R(-1) \\
&&&& \oplus && \oplus \\
&&&&R(-d_2-1) && R(-d_1) \\
0 & \rightarrow & R(-d-1)^2 & \rightarrow & \oplus & \rightarrow & \oplus
& \rightarrow & R & \rightarrow & R/(I+(L)) & \rightarrow & 0 \\
&&&& R(-d_3-1) && R(-d_2) \\
&&&& \oplus && \oplus \\
&&&& R(-d)^2 && R(-d_3)
\end{array}
\]
(The case where $d_1+d_2+d_3$ is odd is analogous.)

\item[b.] If $d_3 > d_1 + d_2 - 2$ then $J = (G_1,G_2)$ is a complete
intersection.  In this case the minimal free resolution of $R/(I+(L))$ is
given by
\[
\begin{array}{ccccccccccccccccccccccccc}
&&&&R(-d_1-1) && R(-1) \\
&&&& \oplus && \oplus \\
0 & \rightarrow & R(-d_1-d_2-1) & \rightarrow & R(-d_2-1) & \rightarrow
& R(-d_1) & \rightarrow  R  \rightarrow  R/(I+(L))  \rightarrow  0\\
&&&& \oplus && \oplus \\
&&&& R(-d_1-d_2) && R(-d_2)
\end{array}
\]
\end{itemize}
\end{corollary}

\begin{proof}
The first part of the corollary is immediate from Lemma \ref{ht2}.

For both (a) and (b) we know that $(I+(L)) = (J+(L))$ where $J$ is
arithmetically Cohen-Macaulay of depth 1.  Hence $R/(I+(L))$ has the same
resolution as
$R/(J+(L))$, either over $R$ or over $R/(L)$.

Let us consider (a).  We know from Corollary \ref{junzomainresult}
that $[I+(L)]/(L) = [J+(L)]/(L) \subset R/(L)$ is an almost complete
intersection, so the same is true of $J \subset R$ since
$\depth R/J = 1$.  Suppose that
$d_1+d_2+d_3$ is even (the case where it is odd is completely analogous).
We have a minimal free resolution (over $R/(L)$) for $R/(J+(L))$ given in
Theorem \ref{mainresult}, so  we thus have a minimal free resolution over
$R$ for $R/J$ given by
\[
\begin{array}{ccccccccccccccc}
&&&& R(-d_1) \\
&&&& \oplus \\
0 & \rightarrow & R(-d)^2 & \rightarrow & R(-d_2) & \rightarrow & R &
\rightarrow & R/J & \rightarrow & 0. \\
&&&& \oplus \\
&&&& R(-d_3)
\end{array}
\]
Then the desired minimal free resolution for $R/(J+(L))$ (and hence
$R/(I+(L))$) is given by the tensor product of this resolution with the
resolution
\[
0 \rightarrow R(-1) \rightarrow R \rightarrow R/(L) \rightarrow 0.
\]
The proof of (b) is trivial
\end{proof}

\begin{remark} It is possible that similar techniques can be used to prove
the Strong Lefschetz property for height three complete intersections (see
Definition \ref{SLPdef}), or to attack either the Weak or Strong Lefschetz
properties for Artinian complete intersections in higher dimensional
rings.  However, a more subtle proof will be needed, as simple examples
show that the {\em degrees} of the syzygies will not be enough to obtain a
contradiction.

Nevertheless, {\em we conjecture that every Artinian complete intersection
in $K[x_0,x_1,x_2]$ has the Strong Lefschetz property}.
\end{remark}

\begin{remark} \label{charp}
What happens in characteristic $p$?  We first note that we cannot expect a
result as strong as the one given in Theorem \ref{mainresult}.  Indeed,
let $A = K[x_1,x_2,x_3]/(x_1^2, x_2^2, x_3^2)$ where $K$ has
characteristic 2.  Let $g = ax_1+bx_2+cx_3$ be a general linear form.
Then $g : A_1 \rightarrow A_2$ is not injective; indeed, $g$ is itself in
the kernel!  A similar observation can be made for $A =
K[x_1,x_2,x_3]/(x_1^4, x_2^4, x_3^4)$, etc.

The main problem here is that the Grauert-M\"ulich theorem does not hold
in characteristic $p$.  There are weaker versions: a Theorem of Ein
(\cite{ein} Theorem 4.1) bounds the splitting type of $\mathcal E$ by a
function of
$c_2({\mathcal E})$.  However, as we saw in the proof of Theorem
\ref{mainresult}, we need the full strength of Grauert-M\"ulich in order
to prove our result.  In the highest degree (at the ``middle'' of the
$h$-vector), the contradiction from the degrees of the syzygies would
not have occurred if this degree had been one greater.  Hence a weaker
version of Grauert-M\"ulich is not good enough with the present techniques.

For example, if $I$ is the complete intersection of three polynomials of
degree 10 in $R$, then one can compute that $c_2({\mathcal E}_{norm}) =
75$, and then Ein's theorem gives that the splitting type is no worse than
$(5,-5)$.  However, that means that the restriction to $\bar R = R/(L)$
has resolution ``no worse'' than
\[
0 \rightarrow \bar R(-10) \oplus \bar R(-20) \rightarrow \bar R(-10)^3
\rightarrow \bar I \rightarrow 0
\]
In particular, it cannot even be excluded that the restriction of $I$ to
$\bar R$ is a complete intersection.  In characteristic zero this is
excluded immediately by our work above (applying the strong
Grauert-M\"ulich theorem) and in fact it follows immediately also from the
main theorem of
\cite{watanabe2}.

\end{remark}

\begin{remark}
\begin{enumerate}
\item The Weak Lefschetz Property says that a general linear form induces
a map of maximal rank on consecutive components.  One might be interested
in a description of the set of (special) linear forms which do {\em not}
give maps of maximal rank.  This is parameterized by the variety of jumping
lines of the bundle $\mathcal E$.

It is interesting to combine the two techniques involved here.  For any
set of distinct lines $\lambda_1, \dots, \lambda_r$ in $\proj{2}$ one can
easily construct bundles having the $\lambda_i$ as jumping lines.  For
example, let $r=3$ and consider complete intersections of type (4,4,4).

On $\lambda_i$, $i = 1,2,3$, choose general points $P_{i,1}, P_{i,2},
P_{i,3}, Q_{i,1}, Q_{i,2}, R_{i,1}, R_{i,2}, R_{i,3}, R_{i,4}$.  Consider
the 4-tuples
\[
(P_{i,1}, P_{i,2}, P_{i,3}, Q_{i,1}), (P_{i,1},
P_{i,2}, P_{i,3}, Q_{i,2}), \hbox{ and } (R_{i,1}, R_{i,2}, R_{i,3},
R_{i,4}).
\]
Choose a general quartic curve $F_1 \in R_4$ containing the 12 points
$(P_{i,1}, P_{i,2}, P_{i,3}, Q_{i,1})$ ($i=1,2,3$).  Choose a general
quartic curve  $F_2 \in R_4$ containing the 12 points $(P_{i,1}, P_{i,2},
P_{i,3}, Q_{i,2})$ ($i=1,2,3$).
Choose a general quartic curve  $F_3 \in R_4$ containing the 12 points
$(R_{i,1}, R_{i,2}, R_{i,3}, R_{i,4})$ ($i=1,2,3$).  (This is possible
since the points were chosen generically.)

Then $(F_1,F_2,F_3)$ is a complete intersection, but its restrictions to
$\lambda_1$, $\lambda_2$ and $\lambda_3$ each have linear syzygies.  Let
${\mathcal E}$ be the bundle constructed from this complete intersection.
Since the restriction to a general line has no smaller than quadratic
syzygies, $\lambda_1, \lambda_2$ and $\lambda_3$ are jumping lines.

\item The bundle $\mathcal E$ used in this section is a Buchsbaum-Rim sheaf.
The interested reader can find a much more extensive treatment of such
sheaves and their properties in \cite{KMNP}, \cite{MP} and \cite{MNP1}.

\end{enumerate}

\end{remark}


\section{Hilbert functions and maximal Betti numbers of algebras with the
Weak Lefschetz property}

In this section we do not  require that $\hbox{char } K = 0$ or that $K$ be
algebraically closed.
We give a complete characterization of the possible Hilbert functions of
algebras with the Weak Lefschetz property.  Furthermore, we show that there
is a sharp upper bound on all of the graded Betti numbers in the minimal
free resolution of an algebra with the Weak Lefschetz property. For the
remainder of this paper we write $R = K[x_0,\dots,x_n]$.

\begin{notation}
If $A=R/I$ is a graded $K$-algebra then we denote the Hilbert function of $A$
by
\[
h_A (t) := \dim_K [R/I]_t.
\]
\end{notation}

The main result of \cite{harima} was to characterize the
Gorenstein sequences (i.e.\ the sequences of integers that can arise as the
Hilbert function of an Artinian Gorenstein ideal) corresponding to
Artinian Gorenstein ideals with the Weak Lefschetz property.  These
turned out to be the so-called Stanley-Iarrobino (SI)-sequences.  As a
consequence, since the height three Gorenstein ideals are well understood
(\cite{buchsbaum-eisenbud}, \cite{diesel} among others), in
$K[x_1,x_2,x_3]$ every Gorenstein sequence occurs as the Hilbert function
of an Artinian ideal with the Weak Lefschetz property.  We now consider the
non-Gorenstein case.

\begin{question} \label{hfq}
Which Hilbert functions (in any
codimension) can occur for ideals whose coordinate rings have the Weak
Lefschetz property?
\end{question}

We will give a complete answer to this question, giving a construction for
an Artinian $K$-algebra with any allowable Hilbert function, having the Weak
Lefschetz property.  Later we will give a bound for the graded Betti numbers
of an Artinian $K$-algebra with the Weak Lefschetz property (Theorem
\ref{thm-bounds}), and we will show that our construction achieves the
bound.  Of course this result includes as a special case the maximal possible
socle type.  However, we have the additional nice result that this maximal
socle type can be read directly from the Hilbert function, so we will
consider the socle type along with the Hilbert function.

Let $A$ be an Artinian graded $K$-algebra  with the Weak Lefschetz
Property,  and let $g$ be a Lefschetz element of $A$.  We make the
following  observations about  the Hilbert function and the socle type of
$A$.

\begin{remark}  \label{hilbnec}
\begin{itemize}
\item[(1)]
Let $d$ be the smallest degree for which  $\times g: A_d\rightarrow
A_{d+1}$ is surjective.  Then the map $\times g: A_j\rightarrow A_{j+1}$ is
also surjective for all $j\geq d$.  This is because we are considering the
natural grading.

\item[(2)]
Hence $\times g: A_j\ra A_{j+1}$ is injective, but not surjective,
for all $j<d$.

\item[(3)]
Let $\underline{h}=(h_0,h_1,\dots,h_s)$ be the Hilbert function of $A$.
From (1) and (2) it follows that
$$
h_0<h_1<\cdots<h_d\geq h_{d+1}\geq\cdots\geq h_s.
$$
In particular, $\underline{h}$ is unimodal and strictly increasing until
it reaches its peak, which is called the {\em Sperner number} of the Hilbert
function of A (\cite{watanabe1}).

\item[(4)]
Thus we see that
there exist integers $u_1,u_2,\dots,u_\ell$ such that
$$
\begin{array}{lll}
h_0<h_1<\cdots<h_{u_1}=\cdots=h_{u_2-1}>h_{u_2}=\cdots=h_{u_3-1} \\
\ \ \ >h_{u_3} \cdots>h_{u_\ell}=\cdots=h_{s}>0.
\end{array}
$$
In particular $u_1=d$.

\item[(5)]
Furthermore from (1) and (2) we have that
the positive part of the first difference of $\underline{h}$, namely
$$
\begin{array}{lllllllllllllllllllllllll}
1, & h_1-h_0, & h_2-h_1, & \cdots, & h_{u_1}-h_{u_1-1},
\end{array}
$$
is the Hilbert function of B=A/(g).  In particular,  this is an
O-sequence.

\item[(6)]
Let $(a_0,\ldots,a_s)$ be the $h$-vector of the socle of $A$. The Hilbert
series  of the socle is called the socle type $S(A,\lambda)$ of $A$,
i.e.\
$$
S(A,\lambda)=\sum_{i=0}^{s}a_i\lambda^i.
$$
We want to compare the socle type with the following polynomial
$$
\Phi_{\underline{h}}(\lambda):=\sum_{i=u_1}^{s}(h_i-h_{i+1})\lambda^i,
$$
where $h_{s+1}=0$.
It can easily be checked from (1), (2) and (4) that
$a_i=0$ for all $i\not\in\{u_2-1,u_3-1,\dots,u_\ell-1,s\}$.
Furthermore we have
$$
a_i\leq h_i-h_{i+1}
$$
for all $i\in\{u_2-1,u_3-1,\dots,u_\ell-1,s\}$.
This follows from
$$
\text{Soc}(A)_i\subset\ker(\times g: A_{i}\ra A_{i+1}),
$$
$\dim\text{Soc}(A)_i=a_i$ and
$\dim\ker(\times g: A_{i}\ra A_{i+1})=h_i-h_{i+1}$.   An Artinian $K$-algebra
for which $a_i = h_i - h_{i-1}$ will be said to have {\em maximal socle
type.}  Notice that the rank of the last free module in the minimal
free resolution of $A$ is equal to $\sum a_i$, the dimension of the socle,
so for an algebra with maximal socle type, this rank is actually equal to
the Sperner number of $A$ (see (3) above).
\end{itemize}
\end{remark}

\medskip

Conditions (3), (4) and (5) give a necessary condition for a Hilbert
function $\underline{h}$ to be the Hilbert function of an Artinian graded
$K$-algebra with the Weak Lefschetz property.  We now show that not only
are these conditions also sufficient, thus characterizing the Hilbert
functions of Artinian $K$-algebras with the Weak Lefschetz property, but
in fact an example exists with the maximal possible socle type, as
described in (6).  We first give the basic construction.

\begin{construction} \label{constr}
Let $\underline{h} = (h_0,h_1,\dots,h_s, h_{s+1} = 0)$ be a finite sequence
of integers satisfying the conditions of (3), (4) and (5) above.  Define
\[
{\bar h}(j) := \max \{ h_j - h_{j-1}, 0 \}.
\]
 Choose Artinian ideals
$$
 {\bar J_1} \subset  {\bar J_2} \subset \ldots  \subset  {\bar J_{\ell}}
\subset {\bar R} :=  K[x_1,\ldots,x_n]
$$
such that $h_{ {\bar R}/ {\bar J_1}} =  {\bar h}$ and $\deg  {\bar J_i} =
h(u_i)$ for all $i = 2,\ldots,{\ell}$.  Now put $J_i = {\bar J_i} R$ for all
$i = 1,\ldots,{\ell}$ and
$$
I := J_1 +  \sum_{i=2}^{\ell} [J_i]_{\geq u_i} + \fm^{s+1},
$$
where $\fm = (x_0,\dots,x_n)$.  Set $A := R/I$.  Note that $J_i$ is not
reduced, but it is the saturated ideal of a zeroscheme $\X_i$.  Furthermore,
we have $\X_1 \supset \X_2 \supset \dots \supset \X_\ell$.
\end{construction}

\begin{proposition} \label{exex}
Let $\underline{h}=(1,h_1,\dots,h_s)$ be a finite sequence of positive
integers.  Then $\underline{h}$ is the Hilbert function of a graded Artinian
$K$-algebra $R/J$ having the Weak Lefschetz property if and only if
$\underline{h}$ is a unimodal O-sequence  such that the positive part of the
first difference is an O-sequence.

Furthermore, let $u_1,\dots,u_\ell$ and $\Phi_{\underline{h}}(\lambda)$ be as
in Remark \ref{hilbnec}.  Then the $K$-algebra $A$ of Construction
\ref{constr} has the Weak Lefschetz property,  Hilbert function
$\underline{h}$ and maximal socle type $\Phi_{\underline{h}}(\lambda)$.
\end{proposition}

\begin{proof}
The necessity is proved in Remark \ref{hilbnec}.  The sufficiency follows
immediately from the claim about Construction \ref{constr}, which we now
prove.

\medskip

\begin{itemize}
\item[(1)]
The Artinian $K$-algebra $A$ has the Weak Lefschetz property:  Let
$B^{(j)}:=R/J_j=\oplus[B^{(j)}]_i$.  We may assume that
$x_0$ is not zero divisor mod $J_j$ for all $j$.
Considering the following commutative diagram
$$
\begin{array}{ccccccc}
[B^{(j)}]_{u_{j+1}-1} & \stackrel{x_0}{\ra} & [B^{(j)}]_{u_{j+1}}
& \ra & [B^{(j+1)}]_{u_{j+1}}\\[2ex]
\parallel & & & & \parallel  \\[2ex]
A_{u_{j+1}-1} & & \stackrel{x_0}{\ra} & & A_{u_{j+1}}
\end{array}
$$
we have, as the proof of Lemma 3.2 in \cite{harima}, that $A$ has the
Weak Lefschetz property.

\item[(2)]
The Hilbert function of $A$ is $\underline{h}$:
First we recall a basic property
of the Hilbert function of a zeroscheme $\Y$ in $\proj{n}$.  Set
$$
\sigma(\Y):=\min\{i\mid \Delta h_{R/I_\Y} (i) = 0\},
$$
where $\Delta h_{R/I_\Y} (i)$ is the first difference of $h_{R/I_\Y} (i)$.
Then it follows that
$$
h_{R/I_\Y} (0)<\cdots<h_{R/I_\Y}
(\sigma(\Y)-1)=h_{R/I_\Y} (\sigma(\Y))=\cdots=\deg \Y,
$$
and we see that
if $\Y^\prime\subset\Y$ then $\sigma(\Y^\prime)\leq\sigma(\Y)$.
Hence from this property we get
$$
h_{B^{(j)}}(i)=h_{u_j}
$$
for all $i\geq u_1$.
Thus since
$A_i=[B^{(1)}]_i$ for all $0\leq i\leq u_2-1$,
$A_i=[B^{(j)}]_i$ for all $u_j\leq i\leq u_{j+1}-1$
and $A_i=(0)$ for all $i\geq s+1$,
we have that the Hilbert function of $A$ coincides with $\underline{h}$.

\item[(3)]
The socle type of $A$ is $\Phi_{\underline{h}}(\lambda)$:
We note that
$$
[\text{Soc}(A)]_{u_{j+1}-1}=[I^{(j+1)}]_{u_{j+1}-1}/[I^{(j)}]_{u_{j+1}-1}.
$$
Furthermore we see that
$$
\begin{array}{lll}
  & \dim\{[I^{(j+1)}]_{u_{j+1}-1}/[I^{(j)}]_{u_{j+1}-1}\} \\[1ex]
= & h_{B^{(j)}}(u_{j+1}-1)-h_{B^{(j+1)}}(u_{j+1}-1) \\[1ex]
= & h_{u_j}-h_{u_{j+1}}.
\end{array}
$$
Thus it follows from part (6) of Remark \ref{hilbnec} that
$$
S(A,\lambda)=\Phi_{\underline{h}}(\lambda).
$$
\end{itemize}
This completes the proof.
\end{proof}

\begin{example} \label{1331}
Not all Artinian ideals in $R$ whose Hilbert functions satisfy the
necessary and sufficient conditions given in Proposition \ref{exex} have
the Weak Lefschetz property.  Indeed, we give a simple example of one which
even has the Hilbert function of a complete intersection but does not have
the Weak Lefschetz property.  We take
\[
I = (x_1^2, x_1x_2, x_1x_3, x_2^3 , x_2^2x_3, x_2x_3^2, x_3^4),
\]
 so $R/I$ has Hilbert function ($1 \ 3
\ 3 \ 1)$.  For any linear form $L$, the element $x_1 \in (R/I)_1$ is in
the kernel of multiplication by $L$, hence the Weak Lefschetz property
fails in passing from degree 1 to degree 2.
\end{example}

A finer invariant of an Artinian $K$-algebra is its minimal
free resolution.  It is probably not possible now to give a set of
necessary and sufficient conditions on the graded Betti numbers for the
existence of an ideal with the Weak Lefschetz property and that set of
Betti numbers.  Even in the Gorenstein case this is open.  However, as in
the Gorenstein case \cite{MN3}, we can give a sharp upper bound for the
graded Betti numbers.  We will do this shortly.

However, we begin with some natural questions, which are the
analogs, for resolutions, of results which we know for Hilbert functions.

\begin{question} \label{wlpq}
\begin{enumerate}
\item Is there a minimal free resolution (meaning only the graded
Betti numbers, not the maps) corresponding to an Artinian ideal with a
Hilbert function allowed by Proposition \ref{exex}, which cannot occur for
an ideal with the Weak Lefschetz property?

\item Are there two Artinian ideals, $I_1$ and $I_2$, which have the same
graded Betti numbers, but one has the Weak Lefschetz property and the
other not?
\end{enumerate}
\end{question}

We answer both of these questions.
First we recall some terminology.

\begin{definition} \label{def-of-LS-BF}
Let $>$ denote the degree-lexicographic order on monomial ideals, i.e.\
$x_1^{a_1}\cdots x_n^{a_n} > x_1^{b_1}\cdots x_n^{b_n}$ if the first
nonzero coordinate of the vector
\[
\left ( \sum_{i=1}^n (a_i - b_i), a_1 - b_1 ,\dots,a_n - b_n \right )
\]
is positive.  Let $J$ be a monomial ideal.  Let $m_1,m_2$ be monomials in
$S$ of the same degree such that $m_1 > m_2$.  Then $J$ is a {\em
lex-segment ideal} if $m_2 \in J$ implies $m_1 \in J$.  When
$\hbox{char}(K) = 0$, we say that $J$ is a {\em Borel-fixed ideal} if
\[
m = x_1^{a_1}\cdots x_n^{a_n} \in J, \ a_i >0, \hbox{ implies }
\frac{x_j}{x_i} \cdot m \in J
\]
for all $1 \leq j < i \leq n$.
\end{definition}

\begin{example}\label{nonoccurring resol}
We first answer the first part of Question \ref{wlpq}.
Let $J \subset K[x_1,x_2,x_3]$ be the lex-segment ideal for the Hilbert
function $(1 \ 3 \ 3 \ 1)$.  Then its minimal free resolution is of the
form
\[
\begin{array}{cccccccccccccccccccccccc}
&& R(-4) && R(-3)^3 && R(-2)^3 \\
&& \oplus && \oplus && \oplus \\
0 & \rightarrow & R(-5)^2 & \rightarrow & R(-4)^5 & \rightarrow & R(-3)^3 &
\rightarrow & J & \rightarrow & 0\\
&& \oplus && \oplus && \oplus \\
&& R(-6) && R(-5)^2 && R(-4)
\end{array}
\]
Now let $I$ be any Artinian ideal in $K[x_1,x_2,x_3]$ with these graded
Betti numbers.  The generators of $I$ in degree 2 have three linear
syzygies.  It is not hard to check (e.g.\ using methods of \cite{BGM}) that
this can only happen if they have a common linear factor (so in particular
there is no regular sequence of length 2 among these three quadrics).  But
then after a change of variables we may assume that this common factor is
$x_1$, and we are in the situation of Example
\ref{1331}.  Hence $R/I$ cannot have the Weak Lefschetz property.  (As an
alternative proof, note that the Socle type is
$\lambda+2\lambda^2+\lambda^3$, so it also follows from Remark
\ref{hilbnec} (6) that it cannot have the Weak Lefschetz property.)
\end{example}

\begin{example}
We now give a (positive) answer to the second part of Question \ref{wlpq}.
H.\ Ikeda  has shown (\cite{ikeda} Example 4.4) that there is a Gorenstein
Artinian $K$-algebra $A = R/I$ with Hilbert function $(1,4,10, 10,4,1)$ and
minimal free resolution
\[
0 \rightarrow {\mathbb F}_4 \rightarrow {\mathbb F}_3 \rightarrow {\mathbb
F}_2 \rightarrow {\mathbb F}_1 \rightarrow R \rightarrow R/I \rightarrow
0,
\]
where
\[
\begin{array}{rcl}
{\mathbb F}_1 & = & R(-3)^{10} \oplus R(-4)^6, \\
{\mathbb F}_2 & = & R(-4)^{15} \oplus R(-5)^{15}, \\
{\mathbb F}_3 & = & R(-5)^6 \oplus R(-6)^{10}, \hbox{ and }\\
{\mathbb F}_4 & = & R(-9).
\end{array}
\]
and not possessing the Weak Lefschetz property.  These graded Betti numbers
are precisely the maximum possible for this Hilbert function among
ideals with the Weak Lefschetz property, and an ideal exists with these
graded Betti numbers and with the Weak Lefschetz property, thanks to
\cite{MN3} Theorem 8.13.
\end{example}

In Example \ref{nonoccurring resol} we saw that the resolution of the
lex-segment ideal (which is known to be extremal among all possible
resolutions with the given Hilbert function \cite{bigatti}, \cite{hulett},
and
\cite{pardue} for $\hbox{char } K > 0$) cannot, in general, be the minimal
free resolution of an ideal with the Weak Lefschetz property, and we gave a
reason for this failure based on the beginning of the resolution, and a
different reason based on the end of the resolution.  This suggests the
following question:

\begin{question} \label{max res question}
Let $\underline{h} = (h_0,h_1,\dots,h_s)$ be a Hilbert function which can
occur for Artinian $K$-algebras with the Weak Lefschetz property (see
Proposition
\ref{exex}).  Is there a maximal possible resolution among Artinian ideals
with the Weak Lefschetz property and Hilbert function $\underline{h}$?
\end{question}

We now answer Question \ref{max res question} by establishing
upper bounds for the graded Betti numbers of an artinian
$K$-algebra with the \WLP\ and exhibiting examples where these
bounds are attained. Note that such bounds were found for artinian
Gorenstein algebras with the \WLP\ in \cite{MN3}. We adapt the
techniques developed there to our problem.

We begin by recalling  \cite{MN3}, Lemma 8.3.

\begin{lemma} \label{lem-tor-seq}
Let $M$ be a graded $R$-module, $\ell \in R$ a linear form.  Then there is an
exact sequence of graded $R$-modules (where $\bar R := R/ \ell R)$:
\[
\begin{array}{c}
\cdots \rightarrow \TT_{i-1}^{\bar R} ((0:_M \ell),K)(-1) \rightarrow
\TT_i^R (M,K) \rightarrow \TT_i^{\bar R} (M/ \ell M ,K) \rightarrow \cdots
\hbox{\hskip 1in} \\
\hfill \cdots \rightarrow \TT_1^R (M,K) \rightarrow \TT_1^{\bar R} (M/ \ell
M,K)  \rightarrow 0.
\end{array}
\]
\end{lemma}

\begin{notation}
Now let $A = R/I$ be an artinian graded $K$-algebra with the \WLP,
and let $g \in [R]_1$ be a Lefschetz element of $A$. Denote by $d$
the end of $A/g A$ and by $a$ the initial degree of $0 : g := 0
:_A g$, i.e. $$ d := \max \{j \in \ZZ \s [A/g A]_j \neq 0 \} $$ $$
a := \min \{j \in \ZZ \s [0 : g]_j \neq 0 \}. $$ Observe that $d
\leq a$. Using the notation of Remark \ref{hilbnec} we have $d= u_1, a =
u_2 -1$.

Moreover, we put ${\bar R} := R/g R$ and define
 $$
  \left [ \tor_i^R (M,K) \right ]_j := \rank \left [
\TT_i^R (M,K) \right ]_j.
$$
\end{notation}

Now we can state the next result.

\begin{proposition} \label{prop-comp-betti}
We have for all integers $i, j$:
\[
\begin{array}{l}
\left [ \tor_i^R (A,K) \right ]_{i+j} = \\ \\

\hskip 1cm
\left \{
\begin{array}{ll}
\left [ \tor_i^{\bar R} (A/g A, K) \right ]_{i+j}
& \hbox{if } j \leq a-2 \\
\leq \left [ \tor_i^{\bar R} (A/g A, K) \right ]_{i+j}
& \hbox{if } j = a-1 \\
\leq \left [ \tor_{i-1}^{\bar R} ( 0 : g, K) \right ]_{i+j-1} +
\left [ \tor_{i}^{\bar R} (A/g A, K) \right ]_{i+j}
& \hbox{if } a \leq j \leq d \\
\leq \left [ \tor_{i-1}^{\bar R} ( 0 : g, K) \right ]_{i+j-1}
& \hbox{if } j = d+1 \\
\left [ \tor_{i-1}^{\bar R} ( 0 : g, K) \right ]_{i+j-1}
& \hbox{if } j \geq d+2
\end{array}
\right.
\end{array}
\]
Furthermore, $\TR_{n+1} (A, K) \cong \TT^{\bar R}_n ( 0: g, K) (-1)$.
\end{proposition}

\begin{proof}
Using $[\TT_i^{\bar R} ((0:_A g),K)]_{i+j} = 0$ if $j < a$ and
$[\TT_i^{\bar R} (A/ gA,K)]_{i+j} = 0$ if $j > d$
the claim follows
by analyzing the exact sequence given in Lemma \ref{lem-tor-seq}.
\end{proof}

Observe that the condition $a \leq j  \leq d$ can only be
satisfied if $a = d$.

Next, we need  an elementary estimate.

\begin{lemma} \label{lem-upper-betti}
Let $M$ be a graded $R$-module. Then we have for all integers $i,
j$:
$$
\left [ \tR_i (M, K) \right ]_{i+j} \leq h_M (j) \cdot \binom{n+1}{i}.
$$
\end{lemma}

\begin{proof}
Put $P := R^{n+1} (-1)$. Then the Koszul complex gives the
following minimal free resolution of $R/\fm \cong K$:
$$
0 \to \wedge^{n+1} P \to \ldots \to \wedge^{i+1} P \to \wedge^{i}
P \to \ldots \to P \to R \to K \to 0.
$$
Thus, $[\TR_i (M, K)]_{i+j}$ is the homology of the complex
$$
\left [ \wedge^{i+1} P \otimes M \right ]_{i+j} \to
\left [ \wedge^{i} P \otimes M \right ]_{i+j} \to
\left [ \wedge^{i-1} P \otimes M \right ]_{i+j}.
$$
Since $\rank [ \wedge^{i} P \otimes M]_{i+j} = h_M (j) \cdot
\binom{n+1}{i}$, the claim follows.
\end{proof}

\begin{notation}
Let $h$ be the Hilbert function of an artinian
$K$-algebra $R/I$.  Then there is a uniquely
determined lex-segment ideal $J \subset R$ such that
$R/J$ has
${h}$ as its Hilbert function.  We define
\[
\beta_{i,j} ({h},R) := \left [ \tor_i^R (R/J,K) \right ]_{i+j}.
\]
\end{notation}

\begin{remark}
The numbers $\beta_{i,j} ({h},R)$ can be computed numerically
without considering lex-segment ideals.  Explicit formulas can be found in
\cite{EK}.
\end{remark}

\begin{theorem}[\cite{bigatti}, \cite{hulett}, \cite{pardue}]
\label{thm-bigatti-hulett}
If  $A = R/I$ is an artinian algebra  then we have for all $i,j \in {\mathbb
Z}$
\[
\left [ \tor_i^R (A,K) \right ]_{i+j} \leq \beta_{i,j} ({h_A},R).
\]
\end{theorem}

In order to construct algebras with the \WLP\ and maximal Betti numbers we
need one more technical result.  In the following lemma, for a graded
module $M$ of finite length we denote by $e(M)$ the last degree in which
$M$ is non-zero.

\begin{lemma} \label{lem-key}
Let ${\bar I} \subset {\bar J} \subset {\bar R} :=
K[x_1,\ldots,x_n]$ be artinian ideals. Put $d := e({\bar R}/{\bar
I})$, $I = {\bar I}R$, $J := {\bar J} R$ and $\fa := I + [J]_{\geq
d + 1}$. Then $\fa + x_0 R = I + x_0 R$ and we have for the graded
Betti numbers of $A := R /\fa$: $$ [\tR_i (A, K)]_j  = \left
\{
\begin{array}{ll}
[\tor_i^{{\bar R}} (A / x_0 A, K)]_j & \mif j \neq i + d \\[2pt]
[\tor_i^{\bar R} (A / x_0 A, K)]_j + k \cdot
\binom{n}{i-1} & \mif j = i + d
\end{array} \right.
$$
where $k := \deg I - \deg J$.
\end{lemma}

\begin{proof}
We proceed in several steps.

(I) Since ${\bar I} \subset {\bar J}$, we get $e({\bar R}/{\bar
J}) \leq e({\bar R}/{\bar I}) = d$. Hence, ${\bar I}$ and ${\bar
J}$ are generated by forms of degree $\leq d+1$. In particular, $[J]_{\geq d
+ 1}$ is generated by forms of degree $d+1$.

The ideals $I + x_0 R$ and $\fa + x_0 R$ differ at most in degrees
$\geq d + 1$. Thus,  the Hilbert functions of $A / x_0 A$ and
${\bar R}/{\bar I}$ agree. It follows that $I + x_0 R = \fa + x_0
R$. In particular, we can write $$ \fa = I + x_0 \cdot
(f_1,\ldots,f_k) $$ where $f_1,\ldots,f_k \in [J]_d$ because $J :
x_0 = J$.

(II) Put $\fb := (f_1,\ldots,f_k) R$, i.e.\ $\fa = I + x_0 \cdot
\fb$. For $j \leq d$, multiplication by $x_0$ factors through two
maps of maximal rank:
$$
\begin{array}{cccccc}
[A]_j & & \stackrel{x_0}{\longrightarrow} & & [A]_{j +1} \\[1pt]
\parallel & & & & \parallel \\[1pt]
[R/I]_j & \stackrel{x_0}{\to} & [R/I]_{j+1} & \to & [R/\fa]_{j+1}.
\end{array}
$$
It follows that
$$
0 :_A x_0 \cong [\fa / I]_d \cong K^k(-d)
$$
and, in particular,  $0 :_A x_0 \cong \Soc A$.

(III) Denote by $g_1,\ldots,g_t$ the minimal generators of $I$.
Let $(r_1,\ldots,r_{t},s_1,\ldots,s_k)^t$ be a syzygy of $\fa$, i.e.\
$$
\sum_{i=1}^t r_i g_i + \sum_{j = 1}^k s_j x_0  f_j = 0.
$$
We can write $r_i = {\bar r_i} + x_0 \tilde{r_i}$ where ${\bar
r_i} \in {\bar R}$ and $\tilde{r_i} \in R$. It follows that
$$
\sum_{i=1}^t {\bar r_i} g_i + x_0 \left [\sum_{i=1}^t \tilde{r_i}  g_i
+  \sum_{j = 1}^k s_j f_j \right ]  = 0.
$$
Comparing coefficients we obtain $\sum_{i=1}^t {\bar r_i} g_i = 0$
and $\sum_{i=1}^t \tilde{r_i}  g_i +  \sum_{i = 1}^k s_j f_j =0$.
Thus, we see that $({\bar r_1},\ldots,{\bar r_t},0,\ldots,0)^t + (x_0
\tilde{r_1},\ldots,x_0 \tilde{r_t},s_1,\ldots,s_k)^t$
is a syzygy of $\fa$ if and only if $({\bar r_1},\ldots,{\bar
r_t})^t$ is a syzygy of $I$ and
$(\tilde{r_1},\ldots,\tilde{r_t},s_1,\ldots,s_k)^t$ is a syzygy of
$I + \fb$.

(IV) Let
$$
0 \to {\bar G_n} \to \ldots {\bar G_2} \stackrel{{\bar
\alpha}}{\longrightarrow} {\bar G_1} \stackrel{{\bar
\beta}}{\longrightarrow} {\bar R} \to {\bar R}/{\bar I} \to 0
$$
be a minimal free resolution of ${\bar R}/{\bar I}$ as ${\bar
R}$-module and let
$$
0 \to F_{n} \to \ldots \to F_1 \to R \to A \to 0
$$
be a minimal free resolution of $A$ as $R$-module.
Tensoring by $\bar R$ gives the complex (with $\bar F_i := F_i
\otimes_R \bar R$)
$$
 0 \rightarrow \bar F_n \rightarrow \cdots
\rightarrow
\bar F_2 \stackrel{\alpha}{\longrightarrow} \bar F_1
\stackrel{\beta}{\longrightarrow}
\bar F_0 \rightarrow R/\fa + x_0 R  \rightarrow 0.
$$
Since $\fa +x_0 R = I + x_0 R$ we get
$$
\ker \beta \cong \ker {\bar \beta} \oplus {\bar R}^k (-d-1).
$$
Step (III) shows that $\im \alpha$ splits as
$$
\im \alpha \cong \im {\bar \alpha} \oplus M \leqno(*)
$$
for some ${\bar R}$-module $M$ such that
$$
\ker \beta / \im \alpha \cong {\bar R}^k (-d-1) / M.
$$
The proof of \cite{MN3}, Lemma 8.3 shows $\ker \beta / \im \alpha
\cong 0 :_A x_0 (-1)$. Using step (II) we obtain the exact sequence
of ${\bar R}$-modules
$$
0 \to M \to {\bar R}^k (-d-1)  \to K^k (- d-1) \to 0.
$$
It implies for all integers $i \geq 0$:
$$
\TT_i^{\bar R} (M, K) \cong K^{k \binom{n}{i+1}}(-d-2-i).
$$
From the proof of \cite{MN3}, Lemma 8.3 we also have for $i \geq
0$:
$$
\TR_{i+2} (A, K) \cong \TT_i^{\bar R} (\im \alpha, K).
$$
Hence, the sequence $(*)$ implies our claim.
\end{proof}

We are now ready for the announced result.

\begin{theorem} \label{thm-bounds} \mbox{ }
\begin{itemize}
\item[(a)] Let $A = R/I$ be a $K$-algebra
with the \WLP\ and denote by ${\bar h}: \ZZ \to \ZZ$ the function defined by
$$
{\bar h}(j) := \max \{ \Delta h_A (j), 0 \}.
$$
Then the graded Betti numbers of $A$ satisfy
$$
\left [ \tR_i (A, K) \right ]_{i+j} \leq
\left \{
\begin{array}{ll}
\beta_{i,j} ({\bar h}, {\bar R}) & \hbox{if $j \leq a-1$} \\
\beta_{i,j} ({\bar h}, {\bar R}) +
\max \{0, - \Delta h_A (j+1)\} \cdot \binom{n}{i-1} &
\hbox{if $a \leq j \leq d$} \\
\max \{0, - \Delta h_A (j+1)\} \cdot \binom{n}{i-1}  & \hbox{if $j \geq d+1$}
\end{array}
\right.
$$
\item[(b)] Let $h: \ZZ \to \ZZ$ be a numerical function such that
there is an artinian algebra $R/J$ having the \WLP\ and $h$ as
Hilbert function. Then there is an artinian algebra $A = R/I$ having
the \WLP\ and $h$ as Hilbert function such that equality is true
in {\rm (a)} for all integers $i, j$.
\end{itemize}
\end{theorem}

\begin{proof}
We first prove (a). Since $g$ is a Lefschetz element of $A$, the Hilbert
function of $A/ g A$
is ${\bar h}$ and the Hilbert function of $0 :_A g$ is given by
$$
h_{0 :_A g} (j) = \max \{ 0, - \Delta h_A (j+1) \}.
$$
Thus, our claim is a consequence of Proposition \ref{prop-comp-betti},
Lemma \ref{lem-upper-betti} and Theorem \ref{thm-bigatti-hulett} (using
\cite{pardue} for the case $\hbox{char } K > 0$).

Now we show (b). We use the notation of Remark \ref{hilbnec}.  Consider the
ideal $I$ of Construction \ref{constr}, and assume furthermore that
\[
[\tor_i^{ {\bar R}}({\bar R}/  {\bar  J_1}, K)]_{i + j} = \beta_{i,j}
({\bar h}, {\bar R}) \ \ \ \hbox{ for all integers $i, j$}.
\]
 Such an ideal ${\bar J_1}$ certainly exists: for example,
we can choose it as a lex-segment ideal.

As in step (I) of the proof of Lemma \ref{lem-key} we see that $I + x_0 R =
J_1 + x_0 R$. An argument as in step (II) of that proof shows that
$$
0 :_A x_0 = \Soc A
$$
and
$$
\rank [0 :_A x_0]_j = \max \{0, - \Delta h (j+1)\}.
$$
It follows that $A$ has the \WLP, $x_0$ is a Lefschetz element for $A$ and
$$
[\tor_i^{ {\bar R}}( 0 :_A x_0 , K)]_{i + j} = \max \{0, - \Delta h (j+1)\}
\cdot \binom{n}{i}.
$$
Moreover, since $A/ x_0 A \cong {\bar R}/ {\bar J_1}$ we have
$$
[\tor_i^{ {\bar R}}(A/x_0 A, K)]_{i + j} = \beta_{i,j} ({\bar h}, {\bar
  R}).
$$
Observe again that $d = u_1$ and $a := u_2 -1 \geq d$. If $a \geq d+1$ all
Betti numbers $[\tR_i (A, K) ]_{i+j}$ are determined by Proposition
\ref{prop-comp-betti} if $j \geq d+2$. Since $[A]_j = [R/J_1]_j$ for
$j \leq a$ we get
$$
[\tR_i (A, K)]_{i+j} = [\tor_i^{ {R}} (R/J_1, K)]_{i+j} =
[\tor_i^{ {\bar R}} (A/x_0 A, K)]_{i+j} \quad \mif j \leq d.
$$
The remaining graded Betti numbers $[\tR_i (A, K)]_{i+d+1}$
can be computed recursively from the Hilbert function of
$A$. (A similar computation can be found on page 4386 of
\cite{N-aBM-divisors}.)

Now let  $a = d$. From the definition of $I$
we immediately obtain
$$
[\tR_i (A, K)]_{i+j} = [tor_i^R(R/(J_1+[J_2]_{\geq a}), K)]_{i+j}  \quad
\mbox{for all} \; j \leq d.
$$
Thus, we know these graded Betti numbers by Lemma \ref{lem-key}.
If $j \geq d+2$ we know $[\tR_i (A, K)]_{i+j}$ by Proposition
\ref{prop-comp-betti}. Thus, the remaining Betti numbers can be computed as
in the previous case.

In any case, we can compute all graded Betti numbers of $A$. The result
shows our claim.
\end{proof}

We would also like to point out that there are Hilbert functions
such that all algebras with that Hilbert function and the \WLP\
have the same (maximal) graded Betti numbers. A similar phenomenon
is true for Gorenstein algebras with the \WLP\ (cf.\ \cite{MN3},
Corollary 8.14).

\begin{corollary}
Let $I \subset R$ be an artinian ideal  such that $A := R/I$
has the \WLP\ and its Hilbert function satisfies
$$
h_A(j) = \binom{n+j}{n} \quad \mbox{for all} \; j \leq d = u_1
\leq u_2 -3
$$
and $u_k + 2 \leq u_{k+1}$ for all $k$ with $2 \leq k < \ell$.
Then the graded Betti numbers of $A$ are
$$
\left [ \tor_i^R (A,K) \right ]_{i+j} = \left \{
\begin{array}{ll}
\binom{n+d}{i+d} & \mif j = d \\[4pt]
- \Delta h_A(u_k)  \cdot \binom{n}{i-1} & \mif j = u_k - 1 \\[2pt]
0 & \mbox{otherwise}.
\end{array} \right.
$$
\end{corollary}

\begin{proof}
By assumption we have $a \geq d+2$. Thus, Lemma \ref{lem-tor-seq}
provides
$$
\left [ \tor_i^R (A,K) \right ]_{i+j} =
\left \{
\begin{array}{ll}
\left [ \tor_i^{\bar R} (A/g A, K) \right ]_{i+j}
& \hbox{if } j \leq d \\
0 & \hbox{if } j = d+1 \\
\left [ \tor_{i-1}^{\bar R} ( 0 : g, K) \right ]_{i+j-1}
& \hbox{if } j \geq d+2
\end{array}
\right.
$$
We may assume that $g = x_0$. Then we get $A/x_0 A \cong {\bar
R}/(x_1,\ldots,x_n)^{d+1}$. Thus, the graded Betti numbers  of
$A/gA$ are known (cf., e.g., the proof of \cite{MN3},
Corollary 8.14). This shows our claim for $j \leq d+1$.

Since $A$ has the \WLP\ we have
$$
\rank [0 :_A x_0]_j = \max \{0, - \Delta h (j+1)\}.
$$
This implies
$$
0 :_A x_0 = \Soc A.
$$
Our claim follows.
\end{proof}


\section{Hilbert functions and maximal Betti numbers of algebras with the
Strong Lefschetz property}

In this section we give some results about a more stringent condition,
namely the Strong Lefschetz property.  Several of our results require that
$\hbox{char } K = 0$, (e.g.\  Proposition \ref{SLP for codim 2}), and we
make this assumption throughout this section.

Not all algebras with the Weak
Lefschetz property possess the Strong Lefschetz property in codimension
$\geq 3$.  We show that nevertheless this {\em does} hold in codimension
two.  Furthermore, we give the surprising result that the {\em same}
characterization of Hilbert functions and maximal graded Betti numbers that
we gave in the last section for the Weak Lefschetz property continues to
hold for the Strong  Lefschetz property.

The conditions for the Hilbert function given in the statement of
Proposition \ref{exex} are automatic in codimension two.  In this case,
interestingly, something much stronger than Proposition \ref{exex} holds.
We first recall the notion of the Strong Lefschetz property.

\begin{definition} \label{SLPdef}
An Artinian ideal $I \subset R$ has the {\em Strong Lefschetz property}
if, for a general linear form $L$ and any $d > 0, i \geq 0$, the map
\[
\times L^d : (R/I)_i \rightarrow (R/I)_{i+d}
\]
has maximal rank.
\end{definition}

Clearly if $R/I$ has the Strong Lefschetz property then it has the Weak
Lefschetz property.  However, there are examples of ideals with the Weak
Lefschetz property which do not have the Strong Lefschetz property.

\begin{example}
We first give a simple example of an ideal with the Weak Lefschetz
property but not the Strong Lefschetz property.  Let $I$ be the lex-segment
ideal with generators
\[
x_1^2,\ x_1x_2,\ x_1x_3^2,\ x_2^3,\ x_2^2x_3^2,\ x_2x_3^3,\ x_3^5.
\]
This has Hilbert function $(1\ 3 \ 4 \ 3 \ 1)$, and one can check that for
multiplication by a general linear form $L$ we have maximal rank between
consecutive components, while $L^2$ has the element $x_1$ in the kernel of
the multiplication from degree 1 to degree 3.
\end{example}

Of much greater interest is the fact that there exist examples of {\em
Gorenstein} ideals with the Weak Lefschetz property but not the Strong
Lefschetz property.  One uses the theory of inverse systems.

\begin{example}
Let $R$ be the ring $K[u,v,x,y,z]$ and let
$f=xu^2+ yuv+zv^2$.  The dual of $f$ gives a Gorenstein algebra with
$h$-vector $(1\ 5 \ 5 \ 1)$ (this can be checked, for instance, with the
computer program Macaulay \cite{macaulay} using the script \
\verb+<l_from_dual+).  This algebra has neither the Weak Lefschetz
property nor the Strong Lefschetz property.

However, now take the polynomial $g = uf$.
It gives an algebra with $h$-vector $(1\ 5\ 6\ 5\ 1)$.
It has the Weak Lefschetz property but {\em not} the Strong Lefschetz
property.

More generally, choose an element $g \in S=[u,v] [f]$ which is, in
particular,  homogeneous in the variables $x,y,z,u,v$.
Let $A$ be the algebra obtained from such a form.  Then for a general
linear form $L$, the map $\times L^{s-2} : A_1 \rightarrow A_{s-1}$
is not bijective.  The key to this goes back to P.\ Gordan and M.\ Noether
\cite{GN}.  They showed that if the Hessian of a form is identically zero
then one of the variables can be eliminated by means of a linear change of
the variables, as long as the number of variable is at most {\em four}.
In dimension 5 or more it is not true, and they gave the above example.
In dimension 5 they claimed that these types of forms are the only cases,
where you have zero Hessian and still all variables are essentially
involved.  Then the fourth author \cite{watanabe3} showed that the zero
Hessian of a form is equivalent to the condition that the map $g^{s-2} :
A_1 \rightarrow A_{s-1}$ does not have full rank.

We believe that in general a polynomial of the above form does give rise to
an Artinian algebra with the Weak Lefschetz property, but have not confirmed
it.
\end{example}

We saw in Example \ref{1331} that for a given Hilbert function in
codimension $\geq 3$ it is possible to find two ideals with that Hilbert
function, one possessing the Weak Lefschetz property and the other not.
In codimension two we have the following, generalizing some results in
\cite{iarrobino}:

\begin{proposition} \label{SLP for codim 2}
Every Artinian ideal in $K[x,y]$ ($\hbox{char } K = 0$) has
the Strong Lefschetz property (and consequently also the Weak Lefschetz
property).
\end{proposition}

\begin{proof}  First suppose that $I$ is a Borel-fixed ideal in
$R=K[x,y]$.
Since $\hbox{char } K = 0$, $I_d$  consists of consecutive monomials from
the first (each $d$).  (Say $x^d$ is the first monomial
and $y^d$ the last.)
So the vector space $R/I_d$ is spanned by the
consecutive monomials from the last.

Let $(h_0, h_1 ,...,  h_s)$ be the Hilbert function of $A=R/I$.
Then it is well known (and easy to see)
that it is unimodal.
Assume first that $h_i \leq h_{i+d}$.
Then
$y^d : (R/I)_i \rightarrow (R/I)_{i+d}$
is injective,
because if a monomial $M$ is in $(R/I)_i$ then
$y^dM$ is in $(R/I)_{i+d}$.
(The point here is that if $M$ is the
$t$-th monomial of $(R/I)_i$
from the last then
$y^dM$ is also the $t$-th monomial of $(R/I)_{i+d}$
from the last. )

Now assume that
 $ h_i \geq h_{i+d}$.
Suppose that a monomial $M$ is in $(R/I)_{i+d}$.
Say $M$ is the $t$-th monomial from the last.
Then the $t$-th monomial of $(R/I)_{i}$ from the last
exists since $h_i > h_{i+d}$.
Let it be $N$.  Then we  have $y^dN=M$.
Thus the map $y^d:(R/I)_i \rightarrow (R/I)_{i+d}$ is
surjective.
Hence we have proved that if $I$ is Borel-fixed in characteristic 0,
then $R/I$ has the Strong Lefschetz property.

In the general case we have the fact that $\hbox{gin}(I)$ is Borel-fixed,
where $\hbox{gin}(I)$ denotes the generic initial ideal of $I$.   It is
easy to see and well known (or see  Proposition 15.12, Eisenbud
\cite{eisenbud}) that
$\hbox{In}(I:y^d) = \hbox{In}(I):y^d$ for $d=1,2,3\dots$, where
$y$ is the last variable with respect to the reverse lexicographic order.
Since the Hilbert function does  not change by passing to $\hbox{gin}(I)$
and  since the Strong Lefschetz property is characterized by the Hilbert
function of $A/(y^d))$, $d=1,2,3\dots$, the general case reduces to the case
of Borel-fixed ideals.
\end{proof}

We have seen that the Strong Lefschetz property is (naturally) a stronger
condition than the Weak Lefschetz property, in the sense that there exist
ideals whose coordinate ring has the Weak Lefschetz property but not the
Strong Lefschetz property.  One would naturally expect that the imposition
of this extra condition would be accompanied by a further restriction on the
possible Hilbert functions (Proposition \ref{exex}) or on the upper bounds on
the graded Betti numbers (Theorem \ref{thm-bounds}).

We now show that any Hilbert function
that occurs for ideals with the Weak Lefschetz property also occurs for
ideals with the Strong Lefschetz property.  The following two results do
not require $\hbox{char } K = 0$.

\begin{proposition} \label{WLP ex has SLP}
Let $K$ be any field.
Let $I$ be the ideal obtained in Construction \ref{constr},
with the further assumption that
$\bar{J_2},\dots,\bar{J_\ell}$ satisfy
\[
h_{\bar{R}/\bar{J_i}}(t) = \Delta h^{(i)}(t)
\]
for all $i=2,\dots,\ell$, where
\[
h^{(i)}(t) := \left \{
\begin{array}{lll}
 \min\{ h_t, h_{u_i} \} & \hbox{ if $t < u_i$,} \\
 h_{u_i}  & \hbox{ otherwise.}
\end{array}
\right.
\]
(Such ideals certainly exist.  For example, we can choose those as
lex-segment ideals.)
Then $A = R/I$ has the Strong Lefschetz property.
\end{proposition}

\begin{proof}
We maintain the notation of Construction \ref{constr} and Proposition
\ref{exex}.  We may assume that  $x_0$ is not a zero divisor mod
$J_j$.   First suppose that $i+d<u_2$.  Then from the proof of Proposition
\ref{exex},  we see that $(A,x_0)$ has the Weak Lefschetz Property.
Hence it follows that the map $\times x_0^d: A_i\ra A_{i+d}$ is injective.

So without loss of generality we may assume that $u_j\leq i+d\leq u_{j+1}-1$
(where $2\leq j\leq\ell$ and $u_{\ell+1}:=s+1)$.
We note that
\[
h_{B^{(j)}}(t) =
\left \{
\begin{array}{ll}
h_t & \hbox{if } \ 0\leq t\leq \sigma(\X_j)-2, \\
h_{u_j} & \hbox{otherwise }.
\end{array}
\right.
\]
Hence we see that
\[
\hbox{the natural map}\
A_i \ra B^{(j)}_i \ \hbox{is} \
\left \{
\begin{array}{ll}
\hbox{bijective} & \hbox{if} \ 0\leq i\leq \sigma(\X_j)-2, \\
\hbox{surjective} & \hbox{if} \ \sigma(\X_j)-1\leq i\leq u_j-1, \\
\hbox{bijective} & \hbox{if} \ u_j\leq i\leq u_{j+1}-1.
\end{array}
\right.
\]
Also we note that
\[
x_0^d: B^{(j)}_i\ra B^{(j)}_{i+d}\ \hbox{is}\
\left \{
\begin{array}{ll}
\hbox{injective} & \hbox{if}\ i\leq\sigma(\X_j)-2, \\
\hbox{bijective} & \hbox{otherwise}.
\end{array}
\right.
\]
Thus, considering the following commutative diagram
$$
\begin{array}{ccc}
A_i & \stackrel{x_0^d}{\longrightarrow} & A_{i+d} \\[2ex]
\downarrow & & \parallel  \\[2ex]
B^{(j)}_i & \stackrel{x_0^d}{\longrightarrow} & B^{(j)}_{i+d}
\end{array}
$$
we have
\[
x_0^d: A_i\ra A_{i+d}\ \hbox{is}\
\left \{
\begin{array}{ll}
\hbox{injective} & \hbox{if}\ i\leq\sigma(\X_j)-2, \\
\hbox{surjective} & \hbox{otherwise}.
\end{array}
\right.
\]
\end{proof}

The next result shows that the bounds on the graded
Betti numbers that were given in Theorem \ref{thm-bounds} are also achieved
by an ideal with the Strong Lefschetz property.

\begin{corollary} \label{SLP hf and betti bounds}
Let $K$ be any field.
A Hilbert function $\underline{h}$ occurs for some graded Artinian
$K$-algebra with the Weak Lefschetz property if and only if it occurs for
one with the Strong Lefschetz property, and these Hilbert functions are
characterized in Proposition \ref{exex}.  Furthermore, the bound on the
graded Betti numbers obtained in Theorem \ref{thm-bounds} is achieved by an
algebra with the Strong Lefschetz property.
\end{corollary}

\begin{proof}
The only thing that needs to be observed is that the extra condition on $\bar
J_1$ imposed in Theorem \ref{thm-bounds}, namely
\[
[\tor_i^{ {\bar R}}({\bar R}/  {\bar  J_1}, K)]_{i + j} = \beta_{i,j}
({\bar h}, {\bar R}) \ \ \ \hbox{ for all integers $i, j$},
\]
can be imposed in the context of Proposition \ref{WLP ex has SLP}: simply
take $\bar J_1$ to be a lex-segment ideal.
\end{proof}

We end with a natural question.

\begin{question}
Is there a set of graded Betti numbers that occurs for algebras with the
Weak Lefschetz property but not the Strong Lefschetz property?
\end{question}

\noindent We conjecture the answer to this question to be ``no."


\end{document}